\RequirePackage[l2tabu,orthodox,abort]{nag}

\documentclass[
11pt,
a4paper,
reqno
]
{amsart}

\usepackage{xcolor}
\usepackage[all,error]{onlyamsmath}

\usepackage[a4paper,left=3cm,right=3cm,top=3cm,bottom=3cm]{geometry}
\usepackage[T1]{fontenc}
\usepackage[utf8]{inputenc}

\usepackage{tikz}
\usepackage{amsmath}
\usepackage{amsfonts}
\usepackage{amssymb}
\usepackage{amsthm}

\theoremstyle{plain}

\newtheorem{theorem}{Theorem}

\theoremstyle{definition}

\usepackage{needspace}

\newcounter{step}
\newcommand{\step}[1][]{%
  \refstepcounter{step}%
  \par\smallskip
  \Needspace{3\baselineskip}
  \noindent\textbf{Step \arabic{step}.}%
  \if\relax\detokenize{#1}\relax\else\ #1\fi
  \par\smallskip\nobreak
  \noindent\ignorespaces
}

\newcommand{\N}{\mathbb{N}}
\newcommand{\R}{\mathbb{R}}
\newcommand{\Z}{\mathbb{Z}}

\newcommand{\C}{\mathbb{C}}

\newcommand{\Lap}{\mathcal{L}}
\newcommand{\Rmax}{R_{\max}}
\newcommand{\Zmax}{Z_{\max}}
\newcommand{\Zdom}{Z_{\mathrm{dom}}}
\newcommand{\mumax}{\mu_{\max}}
\newcommand{\zmin}{z_{\min}}

\newcommand{\e}{\mathrm{e}}

\newcommand{\dd}{\mathop{}\!\mathrm{d}}

\DeclareMathOperator{\re}{Re}
\DeclareMathOperator{\im}{Im}
\DeclareMathOperator{\res}{Res}

\usepackage{mathtools}
\DeclarePairedDelimiter{\lrp}{\lparen}{\rparen}
\DeclarePairedDelimiter{\lrb}{\lbrack}{\rbrack}
\DeclarePairedDelimiter{\lrc}{\lbrace}{\rbrace}
\DeclarePairedDelimiter{\set}{\lbrace}{\rbrace}
\DeclarePairedDelimiter{\abs}{\lvert}{\rvert}
\DeclarePairedDelimiter{\norm}{\lVert}{\rVert}

\usepackage[colorlinks,allcolors=blue]{hyperref}

\title{A counterexample to Abel-type asymptotics for scaled Volterra equations}
\author{Adam Gregosiewicz}
\address{%
  Lublin University of Technology\\
  Nadbystrzycka 38 D\\
  20-618 Lublin\\
  Poland}
\email{a.gregosiewicz@pollub.pl}

\begin{document}

\begin{abstract}
We consider scaled Volterra equations of the form \( f_n + n k*f_n = g \) for
\( n \in \N \), where \( g \) is given and \( f_n \) is sought.
We show that global two-sided Abel-type bounds on a positive kernel \( k \) do
not force the solutions \( f_n \) to converge to zero as \( n \to +\infty \).
More precisely, we construct a continuous strictly positive kernel globally
comparable with the Abel kernel \( x^{-1/2} \), and a continuous strictly
positive \( g \), for which a subsequence of \( (f_n)_{n \in \N} \) diverges to
\( +\infty \) at some point \( x_0 > 0 \).
Consequently, the resolvents associated with the scaled kernels \( nk \) need
not form a generalized approximate identity, in contrast to a couple of
classical results.
\end{abstract}

\maketitle

\section{Introduction}
\label{sec:introduction}

We consider the Volterra integral equation of the second kind
\[
f(x) + k * f(x) = g(x),
\qquad x > 0,
\]
where
\[
k * f(x) \coloneqq \int_0^x k(x-t) f(t) \dd t
\]
whenever the integral converges.
For background and further references we refer to the monographs
\cite{brunner:volterra:2017,gripenberg:volterra:1990}.
Under standard local integrability assumptions on \( k \) and \( g \), the
equation is uniquely solvable on every finite interval.
Its solution can be represented in terms of the resolvent kernel \( r \)
associated with \( k \), that is, the unique function \( r \) satisfying
\[
r + k * r = k.
\]
To wit,
\[
f = g - r * g,
\]
and so the qualitative behaviour of solutions is governed by \( r \).

For each \( n \in \N \), let \( f_n \) denote the solution of the scaled
equation
\[
f_n + nk * f_n = g.
\]
We ask whether \( f_n \to 0 \) as \( n \to +\infty \), at least pointwise and,
perhaps, uniformly on compact subsets of \( (0,+\infty) \).
The degenerate case \( k \equiv 0 \) is irrelevant, since then \( f_n = g \) for
all \( n \).
If \( r_n \) denotes the resolvent kernel of \( nk \), then
\[
f_n = g - r_n * g,
\]
so the problem reduces to the asymptotic behaviour of the family
\( (r_n)_{n \in \N} \).

As \( n \) grows, the memory term \( nk * f_n \) becomes increasingly dominant,
and one may therefore expect the corresponding resolvents to concentrate near
the origin.
This intuition is supported by a formal Laplace-transform computation.
If the Laplace transform \( \Lap k \) of \( k \) converges on \( (0,+\infty) \),
then the resolvent equation for \( r_n \) formally yields
\[
\lrp[\big]{1 + n \Lap k(s)} \Lap r_n(s) = n \Lap k(s),
\qquad s > 0.
\]
In the setting relevant to this paper, where the kernel is positive, we have
\( \Lap k(s) > 0 \) for all \( s > 0 \), and therefore
\[
\Lap r_n(s) = \frac{n \Lap k(s)}{1 + n \Lap k(s)},
\qquad s > 0.
\]
As \( n \to +\infty \), the right-hand side converges to \( 1 \), the Laplace
transform of the Dirac mass at \( 0 \).
Without additional structure, however, this remains only a heuristic.
A standard way to make it rigorous is to show that the resolvents form a
generalized approximate identity: they are nonnegative, their total mass tends
to \( 1 \), and for every \( \delta > 0 \) their mass on \( (\delta,+\infty) \)
tends to zero as \( n \to +\infty \).
Then \( r_n * g \to g \) uniformly on compact subsets of \( (0,+\infty) \) for
every continuous locally bounded \( g \), and consequently \( f_n \to 0 \) in
the same sense.
This leads naturally to the question of which structural assumptions on \( k \)
enforce such behaviour of the resolvent.

Several classical results identify hypotheses under which the resolvent enjoys
positivity or monotonicity properties.
Friedman~\cite{friedman:integral:1963} showed that complete monotonicity of the
kernel is inherited by the resolvent.
For positive kernels that are integrable on \( (0,1) \), continuous and
nonincreasing on \( (0,+\infty) \), and satisfy the quotient monotonicity
condition saying that \( x \mapsto k(x)/k(x+T) \) is nonincreasing on
\( (0,+\infty) \) for every \( T > 0 \), Miller~\cite{miller:volterra:1968}
obtained positivity, continuity, and integrability of the resolvent.
Gripenberg~\cite{gripenberg:positive:1978} sharpened this line of results by
deriving, under additional assumptions on \( k' \), a positive nonincreasing
resolvent.
In the scaled setting, Berrone~\cite{berrone:resolvent:2000} proved that the
same monotonicity hypotheses imply that, for kernels of the form \( nk \), the
family \( (r_n)_{n \in \N} \) is a generalized approximate identity, and hence
\( f_n \to 0 \) uniformly on compact subsets of \( (0,+\infty) \) for every
continuous locally bounded \( g \).
The completely monotone Abel, or Riemann--Liouville, kernels
\( k(x) = x^{-\alpha} \) with \( 0 < \alpha < 1 \) are canonical examples that
satisfy these assumptions.

The Abel kernels suggest also a natural heuristic: for large \( n \), the
decisive feature might be the order of the singularity at \( 0^+ \).
One may therefore expect that if a positive kernel has the same local size near
\( 0^+ \) as \( x^{-\alpha} \), then the corresponding scaled solutions should
again converge to zero.
The purpose of this paper is to show that this expectation is false, already for
\( \alpha = 1/2 \).
In fact, even global two-sided Abel-type bounds do not guarantee the
concentration of the resolvents near the origin.

Our main result reads as follows.

\begin{theorem}
\label{thm:counterexample}
There exist constants \( a, b > 0 \) and continuous strictly positive functions
\( k \) and \( g \) on \( (0,+\infty) \) such that
\[
\frac{a}{\sqrt{x}} \le k(x) \le \frac{b}{\sqrt{x}},
\qquad x > 0,
\]
and the corresponding family \( (f_n)_{n \in \N} \) of unique solutions of
\[
f_n + nk * f_n = g
\]
has a subsequence \( (f_{n_m})_{m \in \N} \) satisfying
\[
\lim_{m \to +\infty} f_{n_m}(x_0) = +\infty
\]
for some \( x_0 > 0 \).
\end{theorem}

Thus neither positivity of the kernel nor two-sided Abel-type bounds, even valid
on the whole half-line, are sufficient to force the convergence of \( f_n \) to
zero.
Consequently, the resolvents associated with the scaled kernels \( nk \) need
not form a generalized approximate identity.
This is noteworthy for two reasons.
First, the kernel \( k \) is globally comparable with the Abel kernel
\( x \mapsto x^{-1/2} \), so the order of the singularity alone does not
determine the asymptotic regime.
Second, the obstruction to convergence comes neither from a loss of positivity
nor from sign changes, but from the absence of additional structural assumptions
of the type appearing in the classical resolvent theory, such as monotonicity.

\subsection{A preliminary example}
\label{sec:preliminary-example}

Before turning to the proof of Theorem~\ref{thm:counterexample}, we record a
much simpler example showing that the convergence of \( (f_n)_{n \in \N} \) to
zero may fail if the kernel is not separated from zero at the origin.

Consider
\[
k(x) = x, \qquad x > 0,
\]
and \( g(x) = 1 \) for \( x > 0 \).
For each \( n \in \N \), let \( f_n \) denote the unique solution of
\[
f_n(x) + n \int_0^x (x-t) f_n(t) \dd t = 1,
\qquad x > 0.
\]
One checks that \( f_n(x) = \cos(\sqrt{n} x) \) for \( x > 0 \).
Indeed, a direct calculation gives
\[
\int_0^x (x-t) \cos(\sqrt{n} t) \dd t = \frac{1 - \cos(\sqrt{n}x)}{n},
\qquad x > 0.
\]
Hence, \( f_n + nk * f_n = 1 \), but \( (f_n)_{n \in \N} \) does not converge
to zero.

This example lies outside the scope of our main result, because the kernel
\( k(x) = x \) vanishes at the origin.
It nevertheless shows that positivity alone does not force the asymptotic
behaviour suggested by the Abel kernels.

\section{Proof of Theorem~\ref{thm:counterexample}}
\label{sec:proof}

Here we prove Theorem~\ref{thm:counterexample}.
We begin by explaining the idea of the construction and then carry out the
argument step by step.

\subsection{Idea of the construction}
\label{sec:construction-idea}

The counterexample is based on a mechanism that becomes visible after passing to
the Laplace transform.
We construct the kernel in the form
\[
k(x) = \frac{c(x)}{\sqrt{x}},
\qquad x > 0,
\]
where \( c \) is bounded above and below by positive constants and is chosen to
be multiplicatively periodic.
The periodicity of \( c \) yields a simple scaling relation for the Laplace
transform of \( k \).
Next, we prove that a certain auxiliary holomorphic function associated with the
Laplace transform has zeros in the right half-plane.
Along a suitable geometric subsequence, these zeros produce poles in the Laplace
transforms of the scaled solutions that move exponentially far to the right.
Using the Bromwich inversion formula, we represent the solution by a contour
integral, shift the contour, and express the solution in terms of residues at
these poles together with a remainder term; the latter is controlled by a
uniform Riemann--Lebesgue argument.
In the resulting residue expansion, the dominant terms come from the poles
corresponding to the rightmost zeros, and among them the zeros of smallest
modulus are made dominant by a suitable choice of \( g \).
Finally, evaluating at a carefully chosen point aligns the oscillatory factors
coming from a distinguished conjugate pair of zeros, which yields divergence.

\subsection{Construction of the counterexample}
\label{sec:construction}

The strategy outlined above is now carried out in twelve steps.
In Steps~\ref{step:c}--\ref{step:zeros}, we construct the kernel and identify
the relevant zeros of the associated auxiliary function.
In Steps~\ref{step:Lapfn}--\ref{step:dominant}, we derive a residue expansion
for a geometric subsequence of solutions and isolate its dominant part.
Finally, in Steps~\ref{step:N}--\ref{step:x0}, we fix the remaining parameters
and choose an evaluation point \( x_0 \) that forces divergence.

\step[Constructing \( c \).]
\label{step:c}%
Let \( u_0 = 2 \) and \( u_1 = 2\pi - 2 \).
Choose \( L \in \N \) so that \( L^2 > u_1 + 1 \), and set
\[
\lambda \coloneqq L^2.
\]
The value of \( L \) will be fixed later.
Accordingly, all objects introduced in this step depend on \( L \), although we
suppress this dependence in the notation.

Since \( 1 < u_0 < u_1 < \lambda - 1 \) and
\[
\cos u_0 = \cos u_1 < 0,
\qquad \sin u_0 > 0,
\qquad \sin u_1 < 0,
\]
we can choose \( \delta > 0 \) small enough that the closed intervals
\[
J_0 \coloneqq [u_0 - \delta, u_0 + \delta],
\qquad J_1 \coloneqq [u_1 - \delta, u_1 + \delta]
\]
are disjoint, contained in \( (1, \lambda) \), and
\begin{equation}
\label{eq:sin_cos_sign}
\sup_{x \in J_0 \cup J_1} \cos x < 0,
\qquad \inf_{x \in J_0} \sin x > 0,
\qquad \sup_{x \in J_1} \sin x < 0.
\end{equation}
We define functions \( \psi_0, \psi_1 \colon [1,\lambda] \to [0,1] \) such that,
for each \( j \in \set{ 0,1 } \), the function \( \psi_j \) is supported in
\( J_j \), linear on \( [u_j-\delta,u_j] \) and on \( [u_j,u_j+\delta] \), and
satisfies
\[
\psi_j(u_j) = 1,
\qquad \psi_j(u_j \pm \delta) = 0.
\]
Explicitly,
\[
\psi_j(x) \coloneqq \max \lrc[\bigg]{ 0, 1 - \frac{\abs{x-u_j}}{\delta} },
\qquad x \in [1,\lambda],\ j \in \set{ 0,1 }.
\]
Each \( \psi_j \) is continuous and vanishes in a neighborhood of \( 1 \) and
\( \lambda \).

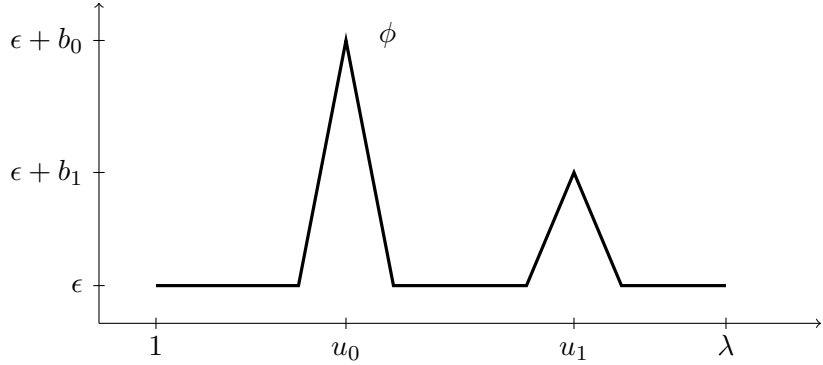
\begin{figure}
\centering
\begin{tikzpicture}[x=2.5cm,y=2.5cm]
  \def\lam{4}
  \def\xzero{2}
  \def\xone{3.2}
  \def\del{0.25}
  \def\eps{0.2}
  \def\bzero{1.3}
  \def\bone{0.6}

  \draw[->] (0.7,0) -- (\lam+0.5,0);
  \draw[->] (0.7,0) -- (0.7,\bzero+0.4);

  \draw (1,0.03) -- (1,-0.03) node[anchor=base,yshift=-10pt] {$1$};
  \draw (\lam,0.03) -- (\lam,-0.03) node[anchor=base,yshift=-10pt] {$\lambda$};
  \draw (\xzero,0.03) -- (\xzero,-0.03) node[anchor=base,yshift=-10pt] {$u_0$};
  \draw (\xone,0.03) -- (\xone,-0.03) node[anchor=base,yshift=-10pt] {$u_1$};

  \draw (0.73,\eps) -- (0.67,\eps) node[left] {$\epsilon$};
  \draw (0.73,{\eps+\bzero}) -- (0.67,{\eps+\bzero}) node[left] {$\epsilon+b_0$};
  \draw (0.73,{\eps+\bone}) -- (0.67,{\eps+\bone}) node[left] {$\epsilon+b_1$};

  \draw[line width=1.2pt]
    (1,\eps) -- (\xzero-\del,\eps) -- (\xzero,{\eps+\bzero}) -- (\xzero+\del,\eps)
    -- (\xone-\del,\eps) -- (\xone,{\eps+\bone}) -- (\xone+\del,\eps) -- (\lam,\eps);
  \node[anchor=west] at ({\xzero+0.12},{\eps+\bzero+0.03}) {$\phi$};
\end{tikzpicture}
\caption{Schematic illustration of \( \phi \).}
\label{fig:phi}
\end{figure}

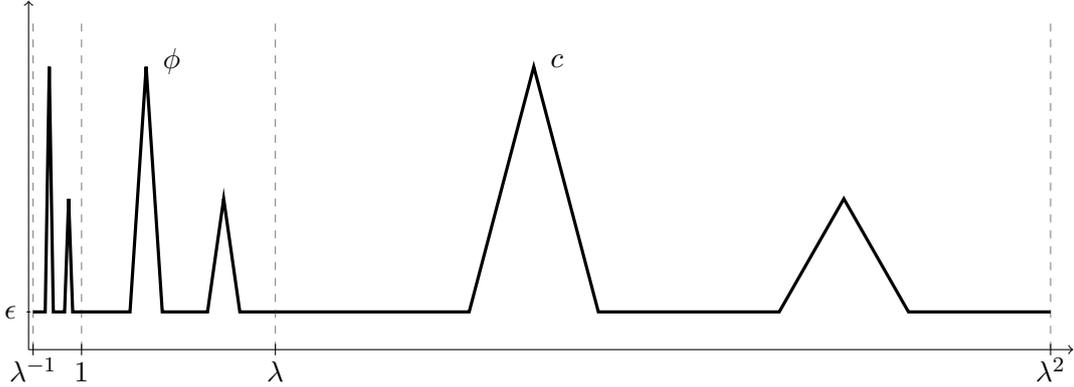
\begin{figure}
\centering
\begin{tikzpicture}[x=0.85cm,y=2.5cm]
  \def\lam{4}
  \def\xzero{2}
  \def\xone{3.2}
  \def\del{0.25}
  \def\eps{0.2}
  \def\bzero{1.3}
  \def\bone{0.6}

  \draw[->] (0.18,0) -- ({\lam*\lam+0.35},0);
  \draw[->] (0.18,0) -- (0.18,{\eps+\bzero+0.35});

  \draw[dashed,gray] (0.25,0) -- (0.25,{\eps+\bzero+0.25});
  \draw[dashed,gray] (1,0) -- (1,{\eps+\bzero+0.25});
  \draw[dashed,gray] (\lam,0) -- (\lam,{\eps+\bzero+0.25});
  \draw[dashed,gray] ({\lam*\lam},0) -- ({\lam*\lam},{\eps+\bzero+0.25});

  \draw (0.25,0.03) -- (0.25,-0.03) node[anchor=base,yshift=-10pt] {$\lambda^{-1}$};
  \draw (1,0.03) -- (1,-0.03) node[anchor=base,yshift=-10pt] {$1$};
  \draw (\lam,0.03) -- (\lam,-0.03) node[anchor=base,yshift=-10pt] {$\lambda$};
  \draw ({\lam*\lam},0.03) -- ({\lam*\lam},-0.03) node[anchor=base,yshift=-10pt] {$\lambda^2$};

  \draw (0.21,\eps) -- (0.15,\eps) node[left] {$\epsilon$};

  \draw[line width=1.2pt]
    (0.25,\eps)
    -- ({\xzero/\lam-\del/\lam},\eps) -- ({\xzero/\lam},{\eps+\bzero}) -- ({\xzero/\lam+\del/\lam},\eps)
    -- ({\xone/\lam-\del/\lam},\eps) -- ({\xone/\lam},{\eps+\bone}) -- ({\xone/\lam+\del/\lam},\eps)
    -- (1,\eps);

  \draw[line width=1.2pt]
    (1,\eps)
    -- (\xzero-\del,\eps) -- (\xzero,{\eps+\bzero}) -- (\xzero+\del,\eps)
    -- (\xone-\del,\eps) -- (\xone,{\eps+\bone}) -- (\xone+\del,\eps)
    -- (\lam,\eps);

  \draw[line width=1.2pt]
    (\lam,\eps)
    -- ({\lam*(\xzero-\del)},\eps) -- ({\lam*\xzero},{\eps+\bzero}) -- ({\lam*(\xzero+\del)},\eps)
    -- ({\lam*(\xone-\del)},\eps) -- ({\lam*\xone},{\eps+\bone}) -- ({\lam*(\xone+\del)},\eps)
    -- ({\lam*\lam},\eps);

  \node[anchor=west] at ({\xzero*\lam+0.1},{\eps+\bzero+0.03}) {$c$};
  \node[anchor=west] at ({\xzero+0.1},{\eps+\bzero+0.03}) {$\phi$};
\end{tikzpicture}
\caption{Function \( c \) as a multiplicatively \(\lambda\)-periodic extension
  of \( \phi \), shown on \( [\lambda^{-1},1] \), \( [1,\lambda] \), and
  \( [\lambda,\lambda^2] \).}
\label{fig:c}
\end{figure}

Fix positive real numbers \( \epsilon \), \( b_0 \), and \( b_1 \), to be
specified later.
Define \( \phi \colon [1,\lambda] \to \R \) by
\[
\phi(x) \coloneqq \epsilon + b_0 \psi_0(x) + b_1 \psi_1(x),
\qquad x \in [1,\lambda].
\]
See Figure~\ref{fig:phi}.
Note that \( \phi(1) = \phi(\lambda) = \epsilon \), since both \( \psi_0 \) and
\( \psi_1 \) vanish near the endpoints of the interval \( [1,\lambda] \).
Finally, we extend \( \phi \) to \( (0,+\infty) \) by multiplicative
\( \lambda \)-periodicity.
That is, for \( x > 0 \), let \( m \in \Z \) be the unique integer such that
\( x/\lambda^m \in [1,\lambda) \), and define
\[
c(x) \coloneqq \phi \lrp[\Big]{ \frac{x}{\lambda^m} }.
\]
For later use, we also extend each \( \psi_j \) to \( (0,+\infty) \) by
multiplicative \( \lambda \)-periodicity; with this convention,
\[
c(x) = \epsilon + b_0 \psi_0(x) + b_1 \psi_1(x), \qquad x > 0.
\]
See Figure~\ref{fig:c} for a schematic illustration.

By construction,
\[
c(\lambda x) = c(x), \qquad x > 0.
\]
Moreover, we prove that \( c \) is continuous on \( (0,+\infty) \).
Indeed, for each \( m \in \Z \), the restriction of \( c \) to
\( (\lambda^m,\lambda^{m+1}) \) is given by \( x \mapsto \phi(x/\lambda^m) \),
and is therefore continuous.
On the other hand, \( c(x) \) converges to \( \phi(\lambda) \) as \( x \)
converges to \( \lambda^m \) from below, and to \( \phi(1) \) as \( x \)
converges to \( \lambda^m \) from above.
Since \( \phi(1) = \phi(\lambda) \), this proves continuity on
\( (0,+\infty) \).
Furthermore, \( c \) is bounded and bounded away from zero, because
\[
\epsilon \le c(x) \le \epsilon + b_0 + b_1, \qquad x > 0.
\]
We define the kernel \( k \) by
\[
k(x) \coloneqq \frac{c(x)}{\sqrt{x}}, \qquad x > 0.
\]

We stress that the specific piecewise linear choice of \( \phi \) is not
essential.
What matters for the argument is only that \( \phi \) has a positive background
level and two independently adjustable bumps supported in the disjoint intervals
\( J_0 \) and \( J_1 \).

\step[Calculating the Laplace transform of \( k \).]
We consider the Laplace transform \( K \) of \( k \), which is given by
\[
K(s) \coloneqq \int_0^{+\infty} \frac{c(x)}{\sqrt{x}} \e^{-sx} \dd x
\]
for every complex number \( s \) in the right half-plane
\( \C_+ = \set{ z \in \C\colon \re z > 0 } \).
Clearly, \( K(s) \) converges for every \( s \in \C_+ \), since
\begin{equation}
\label{eq:Kbound}
\abs{ K(s) }
\le \sup_{x > 0} c(x) \int_0^{+\infty} \frac{1}{\sqrt{x}} \e^{-x\re s} \dd x
\le \frac{\sqrt{\pi}}{\sqrt{\re s}} (\epsilon + b_0 + b_1),
\end{equation}
and, by Morera's theorem, \( K \) is a holomorphic function in \( \C_+ \).
\[
K(\lambda s) = \int_0^{+\infty} \frac{c(x)}{\sqrt{x}} \e^{-\lambda s x} \dd x =
\frac{1}{\sqrt{\lambda}} \int_0^{+\infty} \frac{c(x/\lambda)}{\sqrt{x}} \e^{-s
  x} \dd x.
\]
By the multiplicative periodicity of \( c \), we have \( c(x/\lambda) = c(x) \),
which leads to \( K(\lambda s) = \lambda^{-1/2} K(s) = L^{-1} K(s) \).
Consequently, by induction, we obtain
\begin{equation}
\label{eq:Kmult_periodic}
K(\lambda^m s) = L^{-m} K(s), \qquad m \in \Z,\ s \in \C_+.
\end{equation}

\step[Choosing \( \epsilon \), \( b_0 \) and \( b_1 \) such that
\( K(z_0) = -1 \) for \( z_0 = 1+i \).]
The choice of \( z_0 \) is not essential, but it is convenient because it is
compatible with the earlier choices of \( u_0 \) and \( u_1 \).
Namely, \( \re(\e^{-z_0 x}) = \e^{-x}\cos x \) is negative on both \( J_0 \) and
\( J_1 \), while \( \im(\e^{-z_0 x}) = -\e^{-x}\sin x \) has opposite signs on
these intervals.

By the definition of \( c \), we have
\[
K(z_0) = \epsilon \int_0^{+\infty} x^{-1/2} \e^{-z_0 x} \dd x + b_0 I_0 + b_1
I_1 = \epsilon \sqrt{\pi} z_0^{-1/2} + b_0 I_0 + b_1 I_1,
\]
where \( z_0^{-1/2} \) is defined using the principal branch of the square root,
and
\[
I_j \coloneqq \int_0^{+\infty} \frac{\psi_j(x)}{\sqrt{x}} \e^{-z_0 x} \dd x, \qquad j
\in \set{ 0,1 }.
\]

First we prove that for sufficiently large \( L \) we have
\begin{equation}
\label{eq:quadrants}
\re I_j < 0, \quad j \in \set{ 0,1 }, \qquad \im I_0 < 0, \qquad \im I_1 > 0.
\end{equation}
To this end, note that
\[
I_j = \sum_{m \in \Z} \int_{\lambda^m}^{\lambda^{m+1}}
\frac{\psi_j(x)}{\sqrt{x}} \e^{-z_0 x} \dd x.
\]
In every integral we change variables \( x = \lambda^m y \), which leads to
\[
I_j = \sum_{m \in \Z} \int_1^{\lambda} \frac{\psi_j(\lambda^m
  y)}{\sqrt{\lambda^m y}} \e^{-z_0 \lambda^m y} \lambda^m \dd y = \sum_{m \in
  \Z} L^m \int_1^{\lambda} \frac{\psi_j(y)}{\sqrt{y}} \e^{-z_0 \lambda^m y} \dd
y,
\]
where in the last equality we used the multiplicative \( \lambda \)-periodic
extension of \( \psi_j \).
Denote by \( I_j^0 \) the term of the sum associated to \( m = 0 \).
We have
\[
\re I_j^0 = \int_1^{\lambda} \frac{\psi_j(x)}{\sqrt{x}} \e^{-x} \cos x \dd x,
\qquad \im I_j^0 = - \int_1^{\lambda} \frac{\psi_j(x)}{\sqrt{x}} \e^{-x} \sin x
\dd x.
\]
By Step~\ref{step:c} the support of \( \psi_j \) in \( [1,\lambda] \) equals
\( J_j \).
By~\eqref{eq:sin_cos_sign} it follows that
\[
\re I_j^0 < 0, \quad j \in \set{ 0,1 }, \qquad \im I_0^0 < 0, \qquad \im I_1^0 >
0.
\]
In other words, \eqref{eq:quadrants} holds with \( I_j \) replaced by
\( I_j^0 \).

Now we show that for large \( L \), the integral \( I_j \) is close to
\( I_j^0 \).
We estimate the sum in \( I_j \) independently for \( m < 0 \) and \( m > 0 \).
If \( m \le -1 \), then \( \abs{ \e^{-z_0 \lambda^m y} } = \e^{-\lambda^m y} \le 1
\) for \( y \in [1,\lambda] \), thus for \( j \in \set{ 0,1 } \),
\[
\sum_{m \le -1} L^m \int_1^{\lambda} \frac{\psi_j(y)}{\sqrt{y}} \abs{ \e^{-z_0
    \lambda^m y} } \dd y \le M_j \sum_{m \le -1} L^m = \frac{M_j}{L - 1},
\]
where
\[
M_j = \int_1^{\lambda} \frac{\psi_j(y)}{\sqrt{y}} \dd y, \qquad j \in \set{ 0,1
}.
\]
Analogously, if \( m \ge 1 \), then
\( \abs{ \e^{-z_0 \lambda^m y} } = \e^{-\lambda^m y} \le \e^{-\lambda^m} \) for
\( y \in [1,\lambda] \), and
\[
\sum_{m = 1}^{+\infty} L^m \int_1^{\lambda} \frac{\psi_j(y)}{\sqrt{y}} \abs{
  \e^{-z_0 \lambda^m y} } \dd y \le M_j \sum_{m = 1}^{+\infty} L^m
\e^{-\lambda^m} = M_j \sum_{m=1}^{+\infty} L^m \e^{-L^{2m}}.
\]
Hence, as \( L \to +\infty \), both sums for \( m < 0 \) and \( m > 0 \)
converge to zero (for the second one we can use dominated convergence), and
consequently \( I_j \to I_j^0 \) for \( j \in \set{ 0,1 } \), which
establishes~\eqref{eq:quadrants} for sufficiently large \( L \).
Henceforth, we fix such \( L \).

We are ready to prove the existence of positive \( \epsilon \), \( b_0 \) and
\( b_1 \) such that \( K(z_0) = -1 \).
Let
\[
I_0 = A_0 - iB_0,
\qquad I_1 = A_1 + iB_1
\]
for \( A_0, B_0, A_1, B_1 \in \R \).
By~\eqref{eq:quadrants}, \( A_0 \) and \( A_1 \) are negative, while \( B_0 \)
and \( B_1 \) are positive.
If \( t = B_0 / B_1 \), then \( t > 0 \) and \( I_0 + t I_1 = A_0 + tA_1 \) is
negative.
Denoting \( C = (-A_0 - t A_1)^{-1} \), we obtain \( C I_0 + t C I_1 = -1 \),
that is, \( -1 \) is a linear combination of \( I_0 \) and \( I_1 \) with
positive coefficients.
Moreover, since \( I_0 \) and \( I_1 \) lie in different open quadrants, they
are linearly independent as vectors in \( \C \).
Hence, the linear map \( T\colon \R^2 \to \C \) given by
\[
T(b_0,b_1) \coloneqq b_0 I_0 + b_1 I_1,
\qquad b_0,\ b_1 \in \R
\]
is an isomorphism.
Therefore, the image of the open first quadrant \( (0,+\infty)^2 \) via \( T \)
is an open cone in \( \C \) containing \( -1 \).
Since \( -1 - \epsilon \sqrt{\pi} z_0^{-1/2} \) converges to \( -1 \) as
\( \epsilon \to 0^+ \), it follows that for all sufficiently small
\( \epsilon > 0 \) there exist unique \( b_0, b_1 > 0 \) satisfying
\[
b_0 I_0 + b_1 I_1 = -1 - \epsilon \sqrt{\pi} z_0^{-1/2}.
\]
We fix such \( \epsilon, b_0, b_1 > 0 \) and note that \( K(z_0) = -1 \).

From now on, we set
\[
a = \epsilon,
\qquad b = \epsilon + b_0 + b_1.
\]
Then, by Step~\ref{step:c}, for all \( x > 0 \) we have
\( 0 < a \le c(x) \le b \), or equivalently
\[
0 < a \le \sqrt{x} k(x) \le b.
\]

\step[Proving existence of a rightmost zero of \( H = 1 + K \).]
\label{step:zeros}%
Let \( H\colon \C_+ \to \C \) be defined by
\[
H(s) \coloneqq 1 + K(s), \qquad s \in \C_+.
\]
Denote \( R_0 = 4b^2 \pi \).
Then, for every \( s \) with \( \re s \ge R_0 \), by~\eqref{eq:Kbound} we have
\[
\abs{ K(s) } \le \frac{b \sqrt{\pi}}{\sqrt{R_0}} \le \frac{1}{2},
\]
which implies that
\[
\abs{ H(s) } \ge 1 - \abs{ K(s) } \ge \frac{1}{2}.
\]
This means that \( H \) has no zeros in the half-plane
\( \re s \ge R_0 \).

Next we prove that zeros of \( H \) in the half-plane
\( \re s \ge 1/2 \) cannot have arbitrarily large imaginary part.
Consider a closed interval \( I \) in \( (0,+\infty) \), and let
\( \sigma \mapsto k_{\sigma} \) be the function mapping \( I \) to
\( L^1(0,+\infty) \) given by
\[
k_{\sigma}(x) \coloneqq \frac{c(x)}{\sqrt{x}} \e^{-\sigma x}, \qquad x > 0
\]
for every \( \sigma \in I \).
By dominated convergence, the function is continuous, hence its image
\( \set{ k_{\sigma}\colon \sigma \in I } \) is compact in \( L^1(0,+\infty) \).
Fix \( \eta > 0 \), and choose a finite \( \eta \)-net
\( \set{ k_{\sigma_1}, \dots, k_{\sigma_n} } \) of the image of
\( \sigma \mapsto k_{\sigma} \).
Hence, for every \( \sigma \in I \) there exists \( j \in \set{ 1, \dots, n } \)
such that \( \norm{ k_{\sigma} - k_{\sigma_j} }_{L^1(0,+\infty)} \le \eta \).
Moreover, by the Riemann--Lebesgue lemma applied to each function in the
\( \eta \)-net, there is \( \omega_0 > 0 \) such that
\[
\abs[\bigg]{ \int_0^{+\infty} k_{\sigma_j}(x) \e^{-i \omega x} \dd x } \le
\eta, \qquad \abs{ \omega } \ge \omega_0,\ j \in \set{ 1, \dots, n }.
\]
Since \( \eta \) was arbitrary, this proves that
\[
\int_0^{+\infty} k_{\sigma}(x) \e^{-i \omega x} \dd x
\]
converges to zero as \( \abs{ \omega } \to +\infty \) uniformly in
\( \sigma \in I \).
Thus,
\[
\lim_{\abs{ \omega } \to +\infty} \sup_{\sigma \in I}
  \abs{ K(\sigma + i\omega) }
= \lim_{\abs{ \omega } \to +\infty} \sup_{\sigma \in I}
  \abs[\bigg]{ \int_0^{+\infty} k_{\sigma}(x) \e^{-i \omega x} \dd x }
= 0.
\]
Consequently, for every closed interval \( I \subset (0,+\infty) \) there
exists \( \omega_I > 0 \) such that
\begin{equation}
\label{eq:H-away-from-zero}
\abs{ H(\sigma + i\omega) } \ge 1 - \abs{ K(\sigma + i\omega) } \ge \frac{1}{2}
\end{equation}
for all \( \sigma \in I \) and \( \abs{ \omega } \ge \omega_I \).
Applying this with \( I = [1/2, R_0] \), we see that the set of zeros of \( H \)
in the half-plane \( \re s \ge 1/2 \) is contained in a closed rectangle, and is
therefore finite, as zeros of a nonzero holomorphic function cannot accumulate
inside the domain.

Since \( H(z_0) = 0 \) and \( \re z_0 \ge 1 \), the set of zeros of \( H \) in
the half-plane \( \re s \ge 1/2 \) is nonempty.
Let \( \Rmax \) be the maximum of real parts of all such zeros, that is
\[
\Rmax \coloneqq \max \set{ \re z\colon H(z) = 0,\ \re z \ge 1/2 }.
\]
Then \( \Rmax \ge \re z_0 = 1 \), and the set
\[
\Zmax \coloneqq \set{ z \in \C_+\colon H(z) = 0,\ \re z = \Rmax }
\]
is finite and nonempty.
For \( z \in \Zmax \) denote by \( \mu(z) \) the multiplicity of the zero \( z \),
and let
\[
\mumax \coloneqq \max_{z \in \Zmax} \mu(z), \qquad \Zdom \coloneqq \set{ z \in \Zmax\colon
  \mu(z) = \mumax }.
\]
In other words, \( \mumax \) is the maximum of all multiplicities of zeros of
\( H \) on the line \( \re s = \Rmax \), and \( \Zdom \) is the set of such
zeros with multiplicities equal to \( \mumax \).

We note that there is no zero of \( H \) on the positive real half-line, since
\( K(s) > 0 \) for \( s > 0 \).
Moreover, \( c \) is real-valued, thus \( K(\bar s) = \overline{K(s)} \) and
\( H(\bar s) = \overline{H(s)} \).
In particular, this shows that there is no real number in \( \Zmax \) or
\( \Zdom \), and all zeros in these sets occur in conjugate pairs.

\step[Choosing \( g \) and calculating the Laplace transform of \( f_n \).]
\label{step:Lapfn}%
We define the continuous function \( g \) on \( (0,+\infty) \) by
\[
g(x) \coloneqq x^N, \qquad x > 0,
\]
where \( N \in \N \) will be specified later.
The Laplace transform \( G \) of \( g \) is given by
\[
G(s) = \int_0^{+\infty} x^N \e^{-s x} \dd x = \frac{N!}{s^{N+1}}, \qquad \re s >
0.
\]

For \( n \in \N \), let \( f_n \) be the unique solution to
\( f_n + nk * f_n = g \).
First we prove that the Laplace transform of \( f_n \) exists.
Set \( \gamma = n^{2} R_0 = 4 b^2 \pi n^{2} \), and fix \( x > 0 \).
In the space \( C[0,x] \) define the Bielecki norm
\[
\norm{f}_{\gamma} \coloneqq \sup_{0 \le y \le x} \e^{-\gamma y} \abs{f(y)}, \qquad f \in
C[0,x].
\]
Consider \( g \) as an element of \( C[0,x] \), and let \( T_x \) be the mapping
given by
\[
T_x f(y) \coloneqq g(y) - nk * f(y),
\qquad y \in [0,x],
\]
for all \( f \in C[0,x] \).
If \( f \in C[0,x] \), then \( k * f \) is continuous on \( (0,x] \), by the
dominated convergence theorem, and
\[
\abs{ k * f(y) }
\le b \norm{ f }_{C[0,x]} \int_0^y t^{-1/2} \dd t
= 2 b \norm{ f }_{C[0,x]} \sqrt{y},
\qquad y \in (0,x].
\]
Hence, \( k * f \) extends continuously to \( [0,x] \) with \( k * f(0) = 0 \),
and \( T_x \) maps \( C[0,x] \) into \( C[0,x] \).
Since \( c(y) \le b \) for every \( y > 0 \), we have
\begin{align*}
\e^{-\gamma y} \abs{ T_x f(y) - T_x h(y) }
&\le nb \int_0^y \frac{\e^{-\gamma(y-t)}}{\sqrt{y-t}} \e^{-\gamma t}
   \abs{ f(t) - h(t) } \dd t\\
&\le nb \norm{ f-h }_{\gamma}
   \int_0^{+\infty} \frac{\e^{-\gamma t}}{\sqrt{t}} \dd t\\
&= nb \sqrt{\pi}\, \gamma^{-1/2} \norm{ f-h }_{\gamma}\\
&= 2^{-1} \norm{ f-h }_{\gamma}.
\end{align*}
This proves that the mapping \( T_x \) is a contraction in \( C[0,x] \).
Thus, by the Banach fixed point theorem, there exists a unique fixed point of
\( T_x \) in \( C[0,x] \).
If \( T_x f = f \), then \( f \) satisfies \( f(y) + nk * f(y) = g(y) \) for
every \( y \in (0,x] \).
However, we also have \( f_n(y) + nk * f_n(y) = g(y) \) for such \( y \).
By the uniqueness of solutions to the Volterra equation on \( (0,x] \), it
follows that \( f_n(y) = f(y) \) for every \( y \in (0,x] \), and \( f_n \) has
a continuous extension to \( [0,x] \).
In particular, \( \Phi(x) = \sup_{0 < y \le x} \e^{-\gamma y} \abs{ f_n(y) } \)
is finite.
Consequently, for every \( y \in (0,x] \), by the same argument as above,
\begin{align*}
\e^{-\gamma y} \abs{ f_n(y) }
&\le \e^{-\gamma y} g(y) + n b \int_0^y
\frac{1}{\sqrt{y - t}} \e^{-\gamma(y - t)} \e^{-\gamma t} \abs{ f_n(t) } \dd t\\
&\le M_g + 2^{-1} \Phi(x),
\end{align*}
where \( M_g = \sup_{y > 0} \e^{-\gamma y} g(y) \).
We obtain \( \Phi(x) \le 2 M_g \), which leads to
\( \abs{ f_n(x) } \le 2M_g \e^{\gamma x} \) for every \( x > 0 \).
Therefore, for every \( n \in \N \), the Laplace transform
\[
F_n(s) \coloneqq \int_0^{+\infty} f_n(x) \e^{-s x} \dd x
\]
of \( f_n \) exists for all \( s \in \C_+ \) with \( \re s > \gamma \).

Recall that \( K \) is the Laplace transform of \( k \), hence, taking the
Laplace transform of \( f_n + nk * f_n = g \), we obtain
\( F_n(s)(1 + n K(s)) = G(s) \) for sufficiently large \( \re s \).
Consequently,
\begin{equation}
\label{eq:Fn}
F_n(s) = \frac{G(s)}{1 + n K(s)}
\end{equation}
for all sufficiently large \( \re s \) such that \( 1 + n K(s) \) is nonzero.
Since \( 1 + n K \) is holomorphic in \( \re s > 0 \), this shows that
\( G/(1 + n K) \) is a meromorphic continuation of \( F_n \).
Henceforth, we denote this continuation by \( F_n \).

\step[Choosing a subsequence of \( (F_n)_{n \in \N} \) with poles far to the
right.]
\label{step:subsequence}
We consider the sequence \( (n_m)_{m \in \N} \) given by \( n_m = L^m \) for
\( m \in \N \).
For simplicity, we write \( f_m \) and \( F_m \) instead of \( f_{n_m} \) and
\( F_{n_m} \).
By~\eqref{eq:Kmult_periodic}, we have
\[
1 + n_m K(s) = 1 + L^m K(s) = 1 + K(s/\lambda^m) = H(s/\lambda^m),
\]
and, by~\eqref{eq:Fn},
\begin{equation}
\label{eq:Fm}
F_m(s) = \frac{N!}{s^{N+1} H(s/\lambda^m)}.
\end{equation}
Hence, \( s \in \C_+ \) is a pole of \( F_m \) if and only if
\( z = s/\lambda^m \) is a zero of \( H \).
In particular, the maximal real part of these poles equals
\( \lambda^m \Rmax \), and the poles with this real part are precisely the
points \( s = \lambda^m z \) with \( z \in \Zmax \).

\step[Rescaling the Bromwich integral and shifting the contour.]
\label{step:fm}%
Let
\[
\Gamma = R_0 + 1,
\qquad
\Gamma_m = \lambda^m \Gamma.
\]
By Steps~\ref{step:Lapfn} and~\ref{step:subsequence}, the line
\( \re s = \Gamma_m \) lies to the right of both the abscissa of convergence of
\( F_m \), and all poles of \( F_m \).
Hence, by the Bromwich inversion formula
(see~\cite[Theorem~24.4]{doetsch:laplace:1974}),
\[
f_m(x)
= \frac{1}{2\pi i}
  \int_{\Gamma_m - i\infty}^{\Gamma_m + i\infty} F_m(s) \e^{s x} \dd s,
\qquad x > 0.
\]
Substituting \( s = \lambda^m w \), and using~\eqref{eq:Fm}, we obtain
\begin{equation}
\label{eq:fm_Bromwich_scaled}
f_m(x)
= \frac{N!}{2\pi i} \lambda^{-mN}
  \int_{\Gamma - i\infty}^{\Gamma + i\infty}
  \frac{\e^{\lambda^m w x}}{w^{N+1} H(w)} \dd w,
\qquad x > 0.
\end{equation}

Because the set of zeros of \( H \) in \( \re s \ge 1/2 \) is finite, as proved
in Step~\ref{step:zeros}, and every zero in this half-plane that does not belong
to \( \Zmax \) has real part strictly smaller than \( \Rmax \), we choose
\( \eta \in (0, \Rmax - 1/2] \) such that the set of zeros of \( H \) in
\( \re s \ge \Rmax - \eta \) equals \( \Zmax \).
Let
\[
\sigma = \Rmax - \eta.
\]
Since \( \sigma > 0 \), the integrand in~\eqref{eq:fm_Bromwich_scaled} is
meromorphic in the half-plane \( \re w > 0 \), and its poles in
\( \re w \ge \sigma \) are precisely the elements of \( \Zmax \).

\begin{figure}
\centering
\begin{tikzpicture}[x=1.3cm,y=1.0cm]
  \def\xsig{1.8}
  \def\xstar{3.0}
  \def\xgam{4.4}
  \def\w{2.6}

  \draw[->] (0,-\w-0.7) -- (0,\w+0.7);
  \draw[->] (-0.5,0) -- (5.2,0);

  \draw (0.08,\w) -- (-0.08,\w) node[left] {$\omega$};
  \draw (0.08,-\w) -- (-0.08,-\w) node[left] {$-\omega$};
  \node[anchor=south,below] at (\xsig+0.2,0) {$\sigma$};
  \node[anchor=south,below] at (\xgam+0.2,0) {$\Gamma$};

  \draw[thick] (\xsig,-\w) rectangle (\xgam,\w);
  \node at (4.65,1.7) {$\gamma$};

  \draw[->] (2.55,-\w) -- (3.55,-\w);
  \draw[->] (\xgam,-0.9) -- (\xgam,0.9);
  \draw[->] (3.55,\w) -- (2.55,\w);
  \draw[->] (\xsig,0.9) -- (\xsig,-0.9);

  \draw[densely dashed] (\xstar,-\w-0.7) -- (\xstar,\w+0.5);
  \node[below] at (\xstar,-\w-0.8) {$\re w=\Rmax$};

  \foreach \yy in {0.55,1.05,1.75}{
    \node at (\xstar,\yy) {$\times$};
    \node at (\xstar,-\yy) {$\times$};
  }

  \foreach \xx/\yy in {0.55/3.1,0.85/0.65,1.35/1.35,1.45/2.05}{
    \node at (\xx,\yy) {$\times$};
    \node at (\xx,-\yy) {$\times$};
  }

  \node[below] at (\xsig,-\w) {$\sigma-i\omega$};
  \node[above] at (\xsig,\w) {$\sigma+i\omega$};
  \node[below] at (\xgam,-\w) {$\Gamma-i\omega$};
  \node[above] at (\xgam,\w) {$\Gamma+i\omega$};
\end{tikzpicture}
\caption{Contour \( \gamma \).
  For \( \omega \) large enough, in the half-plane \( \re w \ge \sigma \) all
  zeros (marked by crosses) of \( w \mapsto H(w) \) lie inside the rectangle and
  on the line \( \re w = \Rmax \).
  Outside the rectangle all zeros lie in the strip \( 0 < \re w < \sigma \).}
\label{fig:gamma_m}
\end{figure}
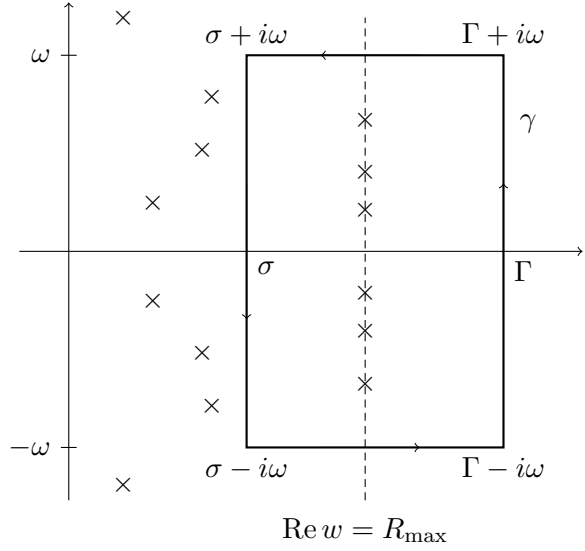

We prove that the integral in~\eqref{eq:fm_Bromwich_scaled} can be rewritten as
the integral over the line \( \re w = \sigma \), plus the sum of the residues of
the integrand at the poles crossed in passing from \( \re w = \Gamma \) to
\( \re w = \sigma \).
To this end, fix \( x > 0 \), let \( \omega > 0 \), and consider the positively
oriented rectangle \( \gamma \) with vertices \( \Gamma \pm i\omega \) and
\( \sigma \pm i\omega \); see Figure~\ref{fig:gamma_m}.
By the residue theorem,
\[
\int_{\gamma} \frac{\e^{\lambda^m w x}}{w^{N+1} H(w)} \dd w
= 2\pi i \sum_{z \in \Zmax} \res_{w = z}
  \lrp[\bigg]{ \frac{\e^{\lambda^m w x}}{w^{N+1} H(w)} }.
\]
If \( \abs{\im w} = \omega \), then \( \abs{w^{-N-1}} \le \omega^{-N-1} \).
Applying~\eqref{eq:H-away-from-zero} with \( I = [\sigma,\Gamma] \), we obtain
\( \omega_0 > 0 \) such that \( \abs{H(u+i\xi)} \ge 1/2 \) for all
\( u \in [\sigma,\Gamma] \) and \( \abs{\xi} \ge \omega_0 \).
Hence, for all \( \omega \ge \omega_0 \),
\[
\abs[\bigg]{
  \int_{\sigma \pm i\omega}^{\Gamma \pm i\omega}
    \frac{\e^{\lambda^m w x}}{w^{N+1} H(w)} \dd w
}
\le \frac{1}{\omega^{N+1}} \e^{\lambda^m \Gamma x}
  \int_{\sigma}^{\Gamma} \frac{\dd u}{\abs{H(u \pm i\omega)}}
\le \frac{2(\Gamma-\sigma)}{\omega^{N+1}} \e^{\lambda^m \Gamma x}.
\]
The last right-hand side tends to \( 0 \) as \( \omega \to +\infty \), which
implies that for all \( x > 0 \) we have
\begin{equation}
\label{eq:fm_residues_reminder}
f_m(x)
= \lambda^{-mN} \sum_{z \in \Zmax} \res_{w = z}
  \lrp[\bigg]{ \frac{N! \e^{\lambda^m w x}}{w^{N+1} H(w)} } + Q_m(x),
\end{equation}
where
\[
Q_m(x)
= \frac{N!}{2\pi i} \lambda^{-mN}
  \int_{\sigma - i\infty}^{\sigma + i\infty}
  \frac{\e^{\lambda^m w x}}{w^{N+1} H(w)} \dd w.
\]

\step[Bounding the remainder term \( Q_m \).]
By the choice of \( \eta \), the function \( \abs{H} \) is positive on the line
\( \re w = \sigma \).
Hence, using continuity and applying~\eqref{eq:H-away-from-zero} with
\( I = \set{\sigma} \), we see that \( \abs{H} \) is bounded away from zero on
that line.
Therefore, denoting
\[
d_{\eta} \coloneqq \inf_{\omega \in \R} \abs{H(\sigma+i\omega)},
\]
we have \( d_{\eta} > 0 \).
Using the definition of \( Q_m \), we get
\[
\abs{Q_m(x)}
\le \frac{N!}{2\pi d_{\eta}} \lambda^{-mN} \e^{\lambda^m \sigma x}
  \int_{-\infty}^{+\infty} \frac{\dd \omega}{(\sigma^2+\omega^2)^{(N+1)/2}},
\qquad x > 0.
\]
Substituting \( \omega = \sigma \tau \) in the last integral, and denoting
\[
\Omega_N \coloneqq \int_{-\infty}^{+\infty} \frac{\dd \tau}{(1+\tau^2)^{(N+1)/2}},
\]
we obtain
\[
\abs{Q_m(x)}
\le \frac{\Omega_N N!}{2\pi d_{\eta}} \cdot
  \frac{\lambda^{-mN} \e^{\lambda^m \sigma x}}{\sigma^N},
\qquad x > 0.
\]
Hence, since \( \sigma = \Rmax - \eta \), there exists \( C > 0 \) such that
\begin{equation}
\label{eq:remainder_bound}
\abs{Q_m(x)} \le C \lambda^{-mN} \e^{\lambda^m (\Rmax - \eta)x},
\qquad m \in \N,\ x > 0.
\end{equation}
The constant \( C \) depends on \( H \), \( \Rmax \), \( \eta \), and \( N \),
but not on \( m \) or \( x \).

\step[Bounding the residues.]
Fix \( z \in \Zmax \) and denote by \( \mu = \mu(z) \) its multiplicity as a zero of
\( H \).
Then there is a holomorphic function \( h \) defined on a neighbourhood of
\( z \) such that \( h(z) \neq 0 \) and \( H(w) = (w-z)^\mu h(w) \) for every
\( w \) in that neighbourhood.
In particular, \( h(z) = H^{(\mu)}(z) / \mu!
\).
Let
\[
a_z(w) \coloneqq \frac{1}{w^{N+1} h(w)}.
\]
Then
\[
\frac{\e^{\lambda^m w x}}{w^{N+1} H(w)}
= \frac{a_z(w)\e^{\lambda^m w x}}{(w-z)^\mu},
\]
so
\begin{align*}
\res_{w = z} \lrp[\bigg]{ \frac{\e^{\lambda^m w x}}{w^{N+1} H(w)} }
&= \frac{1}{(\mu-1)!} \frac{\dd^{\mu-1}}{\dd w^{\mu-1}}
   \lrp[\big]{ a_z(w)\e^{\lambda^m w x} } \Bigr|_{w = z}\\
&= \frac{\e^{\lambda^m z x}}{(\mu-1)!}
   \sum_{j=0}^{\mu-1} \binom{\mu-1}{j} a_z^{(j)}(z) (\lambda^m x)^{\mu-1-j},
\end{align*}
for every \( x > 0 \).
Since \( a_z^{(j)}(z) \) is independent of \( m \), for
\( j \in \set{ 1, \dots, \mu-1 } \) we have
\[
a_z^{(j)}(z) (\lambda^m x)^{\mu-1-j}
= (\lambda^m x)^{\mu-1} O(\lambda^{-m})
\]
as \( m \to +\infty \).
Thus, for every \( x > 0 \) we have
\[
\res_{w = z} \lrp[\bigg]{ \frac{N! \e^{\lambda^m w x}}{w^{N+1} H(w)} }
= \lambda^{m(\mu-1)} \e^{\lambda^m z x} \lrp[\big]{ B_z(x) + O(\lambda^{-m}) }
\]
as \( m \to +\infty \), where
\[
B_z(x)
\coloneqq \frac{N! a_z(z)}{(\mu-1)!} x^{\mu-1}
= \frac{N!\mu}{z^{N+1} H^{(\mu)}(z)} x^{\mu-1}.
\]
Since \( a_z(z) \neq 0 \), we have \( B_z(x) \neq 0 \) for every \( x > 0 \).
Consequently,
\begin{equation}
\label{eq:residue_bound}
\lambda^{-mN} \res_{w = z}
  \lrp[\bigg]{ \frac{N! \e^{\lambda^m w x}}{w^{N+1} H(w)} }
= \lambda^{-m(N+1-\mu)} \e^{\lambda^m z x} \lrp[\big]{ B_z(x) + O(\lambda^{-m}) }.
\end{equation}

\step[Extracting the dominant term in \( f_m \).]
\label{step:dominant}%
By~\eqref{eq:fm_residues_reminder} and~\eqref{eq:residue_bound}, for every
\( x > 0 \) we have
\begin{align*}
f_m(x)
&= \sum_{z \in \Zmax} \lambda^{-m(N+1-\mu(z))} \e^{\lambda^m z x}
   \lrp[\big]{ B_z(x) + O(\lambda^{-m}) } + Q_m(x)\\
&= \lambda^{-m(N+1-\mumax)} \e^{\lambda^m \Rmax x}
   \lrb[\Bigg]{ \sum_{z \in \Zdom} \e^{i\lambda^m \im(z) x} B_z(x) + o(1) }
   + Q_m(x)
\end{align*}
as \( m \to +\infty \).
Indeed, after factoring out \( \lambda^{-m(N+1-\mumax)} \e^{\lambda^m \Rmax x} \),
the terms with \( z \in \Zmax \setminus \Zdom \) tend to zero because
\( \mu(z) < \mumax \), and for \( z \in \Zdom \) the error term
\( O(\lambda^{-m}) \) also tends to zero; since \( \Zmax \) is finite, all
these contributions are absorbed into the \( o(1) \) term.
Combining this with~\eqref{eq:remainder_bound}, we get
\[
f_m(x)
= \lambda^{-m(N+1-\mumax)} \e^{\lambda^m \Rmax x} \lrb[\big]{ S_m(x) + E_m(x) },
\]
where
\[
S_m(x) \coloneqq \sum_{z \in \Zdom} \e^{i\lambda^m \im(z) x} B_z(x),
\]
and \( E_m(x) \) converges to zero as \( m \to +\infty \) for every \( x > 0 \).
The contribution of \( Q_m(x) \) is absorbed into \( E_m(x) \), since
\( \mumax \ge 1 \) and \( \eta > 0 \).
Note that \( S_m \) is real-valued, since \( \Zdom \) consists of pairs of
conjugate zeros, and \( B_{\bar z}(x) = \overline{B_z(x)} \).
Thus, since \( f_m \) and \( S_m \) are real-valued, so is \( E_m \).

\step[Choosing \( N \) so that the terms corresponding to zeros with minimal
modulus dominate \( S_m \).]
\label{step:N}%
Let \( \zmin \) be the unique zero in \( \Zdom \) with positive imaginary part
and minimal modulus.
For every \( z \in \Zdom \setminus \set{ \zmin, \overline{\zmin} } \) we have
\[
\frac{\abs{ B_z(x) }}{\abs{ B_{\zmin}(x) }} = \frac{\abs{ H^{(\mumax)}(\zmin)
  }}{\abs{ H^{(\mumax)}(z) }} \lrp[\bigg]{ \frac{\abs{ \zmin }}{\abs{ z }}
}^{N+1}, \qquad x > 0.
\]
By the choice of \( \zmin \), the right-hand side converges to zero as
\( N \to +\infty \).
Hence, for sufficiently large \( N \), we have
\begin{equation}
\label{eq:sumBz}
\sum_{z \in \Zdom \setminus \set{ \zmin, \overline{\zmin} }} \abs{ B_z(x) } \le
\frac{1}{2} \abs{ B_{\zmin}(x) }, \qquad x > 0.
\end{equation}

Let \( \theta = \arg \zmin \), so that \( \zmin = \abs{\zmin} \e^{i \theta} \),
and write \( B_{\zmin}(1) = \abs{B_{\zmin}(1)} \e^{i \theta_N} \).
Since the argument of \( B_{\zmin}(1) \) depends on \( N \) only through the
factor \( (\zmin)^{-(N+1)} \), there exists \( \tau \in \R \), independent of
\( N \), such that \( \theta_N \) equals \( \tau - (N+1)\theta \) modulo
\( 2\pi \).
If \( \theta/\pi \) is irrational, then the sequence \( (\theta_N)_{N\in\N} \)
is equidistributed modulo \( 2\pi \).
Hence, infinitely many terms lie within distance at most \( \pi/3 \) from
\( 0 \).
On the other hand, if \( \theta/\pi \) is rational, then
\( (\theta_N)_{N\in\N} \) is periodic modulo \( 2\pi \), and therefore takes
finitely many equally spaced values on the circle \( \R / 2\pi\Z \).
Since \( \zmin \) is not a real number, \( \theta \) is not an integer multiple
of \( \pi \), and consequently the set of values of the sequence contains at
least three distinct points.
Therefore, the spacing between consecutive points is at most \( 2\pi/3 \), which
implies that at least one of these points lies within distance at most
\( \pi/3 \) from \( 0 \).
Thus, in either case, there exist infinitely many \( N \) satisfying
\[
\cos \theta_N \ge \frac{1}{2}.
\]

Fix such an \( N \) so large that~\eqref{eq:sumBz} holds.
Then
\[
\abs{ B_{\zmin}(x) } = \abs{ B_{\zmin}(1) x^{\mumax-1} } \le 2 \abs{ \re
  B_{\zmin}(1) x^{\mumax-1} } = 2 \re B_{\zmin}(x), \qquad x > 0,
\]
and
\begin{equation}
\label{eq:sumBzRe}
\sum_{z \in \Zdom \setminus \set{ \zmin, \overline{\zmin} }} \abs{ B_z(x) } \le
\re B_{\zmin}(x), \qquad x > 0.
\end{equation}

\step[Choosing \( x_0 \) such that \( f_m(x_0) \) diverges.]
\label{step:x0}%
Let
\[
x_0 \coloneqq \frac{2\pi}{\im \zmin}.
\]
Since \( \im \zmin > 0 \), we have \( x_0 > 0 \).
Moreover, \( \e^{i \lambda^m \im(\zmin) x_0} = \e^{2\pi i \lambda^m} = 1 \) for
all \( m \in \N \).
Hence,
\begin{align*}
\e^{i\lambda^m \im(\zmin) x_0}B_{\zmin}(x_0)
  + \e^{i\lambda^m \im(\overline{\zmin}) x_0} B_{\overline{\zmin}}(x_0)
&= B_{\zmin}(x_0) + \overline{B_{\zmin}(x_0)}\\
&= 2 \re B_{\zmin}(x_0),
\end{align*}
and by~\eqref{eq:sumBzRe}, it follows that
\[
S_m(x_0) \ge 2 \re B_{\zmin}(x_0) - \sum_{z \in \Zdom \setminus \set{ \zmin,
    \overline{\zmin} }} \abs{ B_z(x_0) } \ge \re B_{\zmin}(x_0).
\]
Since \( E_m(x_0) \) converges to zero as \( m \to +\infty \), we have
\[
S_m(x_0)+E_m(x_0)\ge S_m(x_0) - \abs{ E_m(x_0) } \ge \frac{1}{2} \re
B_{\zmin}(x_0) > 0
\]
for all sufficiently large \( m \).
For these \( m \), we obtain
\[
f_m(x_0) \ge \frac{1}{2} \re B_{\zmin}(x_0) \lambda^{-m(N+1-\mumax)}
\e^{\lambda^m \Rmax x_0}.
\]
Here \( \re B_{\zmin}(x_0) > 0 \), \( \Rmax > 0 \), and \( \lambda > 1 \), which
shows that the right-hand side tends to \( +\infty \), and completes the proof
of Theorem~\ref{thm:counterexample}.

\bibliographystyle{plain}
\bibliography{references.bib}

\end{document}